\documentclass{gtart}

\def\ifplaintex{\expandafter\ifx\csname documentclass\endcsname\relax}

\ifplaintex 
\else
\fi

\expandafter\ifx\csname beginpicture\endcsname\relax
\expandafter\ifx\csname documentclass\endcsname\relax
\input pictex \else
\input prepictex \input pictex \input postpictex \fi\fi

\def\gt{{\mathsurround=0pt\it $\cal G\mskip-2mu$eometry \&\ 
$\cal T\!\!$opology}}        

\def\gtp{{\mathsurround=0pt\it $\cal G\mskip-2mu$eometry \&\ 
$\cal T\!\!$opology $\cal P\!$ublications}}  

\def\lognumber#1{\def\thelognumber{#1}}
\def\volumenumber#1{\def\thevolumenumber{#1}}
\def\papernumber#1{\def\thepapernumber{#1}}
\def\volumeyear#1{\def\thevolumeyear{#1}}

\def\pagenumbers#1#2{\def\startpage{#1}\def\finishpage{#2}}
\def\published#1{\def\publishdate{#1}}
\def\proposed#1{\def\theproposer{#1}}
\def\seconded#1{\def\theseconders{#1}}
\def\received#1{\def\receiveddate{#1}}
\def\revised#1{\def\reviseddate{#1}}
\def\accepted#1{\def\accepteddate{#1}}

\long\def\asciiabstract#1{\long\def\theasciiabstract{#1}}

\let\\\par\let\thelognumber\relax
\let\thevolumenumber\relax\let\thepapernumber\relax
\let\thevolumeyear\relax\let\thesamplenumber\relax\let\startpage\relax
\let\finishpage\relax\let\publishdate\relax\let\receiveddate\relax
\let\reviseddate\relax\let\accepteddate\relax\let\theasciititle\relax
\let\theasciiauthors\relax
\let\theasciiabstract\relax
\let\theasciiemail\relax\let\theshortauthors\relax\let\theshorttitle\relax

\long\def\maketitlep{   

\count0=\startpage

\gt\hfill      
\beginpicture
\setcoordinatesystem units <0.33truein, 0.33truein> point at 2.2 0.9
\setplotsymbol ({$\cal G$})
\plotsymbolspacing=9truept
\circulararc 315 degrees from 0 1 center at 0 0
\setplotsymbol ({$\cal T$})
\circulararc 315 degrees from 1 -1 center at 1 0
\endpicture
\break
{\small\ifx\thesamplenumber\relax 
Volume \else Sample
\fi\thevolumenumber\ (\thevolumeyear)
\startpage--\finishpage\nl
Published: \publishdate}
\vglue 0.5truein plus 0.4fil minus 0.1truein

{\parskip=0pt\leftskip 0pt plus 1fil\def\\{\par\smallskip}{\ifplaintex\large
\else\Large\fi\bf\thetitle}\par\medskip}   

\vglue 0pt plus 0.1fil 

{\parskip=0pt\leftskip 0pt plus 1fil\def\\{\par}{\sc\theauthors}
\par\medskip}

\vglue 0pt plus 0.1fil 

{\small\parskip=0pt\let\newline\\
{\leftskip 0pt plus 1fil\def\\{\par}{\sl\theaddress}\par}
\expandafter\ifx\theemail\relax    
\relax\else\vglue 5pt plus 0.02fil minus 2pt\def\\{\stdspace{\rm 
and}\stdspace} 
\cl{Email:\stdspace\tt\theemail}\fi
\ifx\theurl\relax                  
\relax\else\vglue 5pt plus 0.02fil minus 2pt\def\\{\stdspace{\rm 
and}\stdspace}
\cl{URL:\stdspace\tt\theurl}\fi\par}

\vglue 7pt plus 0.3fil minus 3pt

{\bf Abstract}
\vglue 5pt plus 0.1fil minus 2pt

\theabstract

\vglue 7pt plus 0.3fil minus 3pt

{\bf AMS Classification numbers}\quad Primary:\quad \theprimaryclass

Secondary:\quad \thesecondaryclass

\vglue 5pt plus 0.3fil minus 2pt

{\bf Keywords}\quad \thekeywords

\vglue 10pt plus 0.5fil minus 5pt

{\small  Proposed: \theproposer\hfill Received: \receiveddate\nl
Seconded: \theseconders\hfill 
\ifx\reviseddate\relax                         
Accepted: \accepteddate                        
\else
Revised: \reviseddate                          
\fi}
\eject
}       

\let\maketitlepage\maketitlep
\let\maketitle\maketitlepage

\font\phead=cmsl9 scaled 950
\font=cmsl9 scaled 1050
\font\pnum=cmbx10 scaled 913
\font\lnum=cmbx10 
\font\pfoot=cmsl9 scaled 950
\font=cmsl9 scaled 1050
\ifplaintex
\headline{\vbox to 0pt{\vskip -4.5mm\line{\small\phead\ifnum
\count0=\startpage ISSN 1364-0380 (on line)
1465-3060 (printed) \hfill {\pnum\folio}\else\ifodd\count0\def\\{ }%
\ifx\theshorttitle\relax\thetitle\else\theshorttitle\fi\hfill{\pnum\folio}
\else\def\\{ and }{\pnum\folio}\hfill\ifx\theshortauthors\relax\theauthors
\else\theshortauthors\fi\fi\fi}\vss}}
\footline{\vbox to 0pt{\vglue 0mm\line{\small\pfoot\ifnum\count0=\startpage
\copyright\ \gtp\hfill\else
\gt, Volume \thevolumenumber\ (\thevolumeyear)\hfill\fi}\vss
}}
\else
\makeatletter
\def\@oddhead{{\small\lhead\ifnum\count0=\startpage ISSN 1364-0380 (on line)
1465-3060 (printed) \hfill {\lnum\number\count0}\else\ifodd\count0
\def\\{ }\ifx\theshorttitle\relax \thetitle \else\theshorttitle\fi\hfill
{\lnum\number\count0}\else\def\\{ and }{\lnum\number\count0}
\hfill\ifx\theshortauthors\relax 
\theauthors\else\theshortauthors\fi\fi\fi}}\def\@evenhead{\@oddhead}
\def\@oddfoot{\small\lfoot\ifnum\count0=\startpage\copyright\ \gtp\hfill\else
\gt, Volume \thevolumenumber\ (\thevolumeyear)\hfill\fi}
\def\@evenfoot{\@oddfoot}
\makeatother
\fi

\newwrite\gtoutfile
\long\gdef\makeheadfile{  
{\def\\{, }\def\s{ }
\immediate\openout\gtoutfile head.xxx
\immediate\write\gtoutfile{To: math@arxiv.org}
\immediate\write\gtoutfile{Subject: put or rep NNNNN:pppp}
\immediate\write\gtoutfile{--text follows this line--}
\immediate\write\gtoutfile{Proxy-for: \ifx\theasciiauthors\relax
\theauthors\else\theasciiauthors\fi\s<\ifx\theasciiemail\relax\theemail\else\theasciiemail\fi>}
\immediate\write\gtoutfile{\noexpand\\}
\immediate\write\gtoutfile{Authors: \ifx\theasciiauthors\relax
\theauthors\else\theasciiauthors\fi}
{\def\\{ }\immediate\write\gtoutfile{Title: \ifx\theasciititle\relax
\thetitle\else\theasciititle\fi}}
\immediate\write\gtoutfile{Subj-class: GT or SG or MG etc}
\immediate\write\gtoutfile{MSC-class: \theprimaryclass\ifx\thesecondaryclass\relax\else, \thesecondaryclass\fi}
\immediate\write\gtoutfile{Journal-ref: Geom. Topol. \thevolumenumber
(\thevolumeyear) \startpage-\finishpage}
\immediate\write\gtoutfile{Comments: Published by Geometry and Topology at}
\immediate\write\gtoutfile{\s\s http://www.maths.warwick.ac.uk/gt/GTVol\thevolumenumber/paper\thepapernumber.abs.html}
\immediate\write\gtoutfile{\noexpand\\}
\immediate\write\gtoutfile{}
\ifx\theasciiabstract\relax
\immediate\write\gtoutfile{\theabstract}\else
\immediate\write\gtoutfile{\theasciiabstract}\fi
\immediate\write\gtoutfile{}
\immediate\write\gtoutfile{\noexpand\\}
\immediate\write\gtoutfile{}
\immediate\closeout\gtoutfile}}  

\def\maketitlepage{\maketitlep\makeheadfile}
\let\maketitle\maketitlepage

\def\ifplaintex{\expandafter\ifx\csname documentclass\endcsname\relax}

\ifplaintex 
\else
\fi

\expandafter\ifx\csname beginpicture\endcsname\relax
\expandafter\ifx\csname documentclass\endcsname\relax
\input pictex \else
\input prepictex \input pictex \input postpictex \fi\fi

\def\gt{{\mathsurround=0pt\it $\cal G\mskip-2mu$eometry \&\ 
$\cal T\!\!$opology}}        

\def\gtp{{\mathsurround=0pt\it $\cal G\mskip-2mu$eometry \&\ 
$\cal T\!\!$opology $\cal P\!$ublications}}  

\def\lognumber#1{\def\thelognumber{#1}}
\def\volumenumber#1{\def\thevolumenumber{#1}}
\def\papernumber#1{\def\thepapernumber{#1}}
\def\volumeyear#1{\def\thevolumeyear{#1}}

\def\pagenumbers#1#2{\def\startpage{#1}\def\finishpage{#2}}
\def\published#1{\def\publishdate{#1}}
\def\proposed#1{\def\theproposer{#1}}
\def\seconded#1{\def\theseconders{#1}}
\def\received#1{\def\receiveddate{#1}}
\def\revised#1{\def\reviseddate{#1}}
\def\accepted#1{\def\accepteddate{#1}}

\long\def\asciiabstract#1{\long\def\theasciiabstract{#1}}

\let\\\par\let\thelognumber\relax
\let\thevolumenumber\relax\let\thepapernumber\relax
\let\thevolumeyear\relax\let\thesamplenumber\relax\let\startpage\relax
\let\finishpage\relax\let\publishdate\relax\let\receiveddate\relax
\let\reviseddate\relax\let\accepteddate\relax\let\theasciititle\relax
\let\theasciiauthors\relax
\let\theasciiabstract\relax
\let\theasciiemail\relax\let\theshortauthors\relax\let\theshorttitle\relax

\long\def\maketitlep{   

\count0=\startpage

\gt\hfill      
\beginpicture
\setcoordinatesystem units <0.33truein, 0.33truein> point at 2.2 0.9
\setplotsymbol ({$\cal G$})
\plotsymbolspacing=9truept
\circulararc 315 degrees from 0 1 center at 0 0
\setplotsymbol ({$\cal T$})
\circulararc 315 degrees from 1 -1 center at 1 0
\endpicture
\break
{\small\ifx\thesamplenumber\relax 
Volume \else Sample
\fi\thevolumenumber\ (\thevolumeyear)
\startpage--\finishpage\nl
Published: \publishdate}
\vglue 0.5truein plus 0.4fil minus 0.1truein

{\parskip=0pt\leftskip 0pt plus 1fil\def\\{\par\smallskip}{\ifplaintex\large
\else\Large\fi\bf\thetitle}\par\medskip}   

\vglue 0pt plus 0.1fil 

{\parskip=0pt\leftskip 0pt plus 1fil\def\\{\par}{\sc\theauthors}
\par\medskip}

\vglue 0pt plus 0.1fil 

{\small\parskip=0pt\let\newline\\
{\leftskip 0pt plus 1fil\def\\{\par}{\sl\theaddress}\par}
\expandafter\ifx\theemail\relax    
\relax\else\vglue 5pt plus 0.02fil minus 2pt\def\\{\stdspace{\rm 
and}\stdspace} 
\cl{Email:\stdspace\tt\theemail}\fi
\ifx\theurl\relax                  
\relax\else\vglue 5pt plus 0.02fil minus 2pt\def\\{\stdspace{\rm 
and}\stdspace}
\cl{URL:\stdspace\tt\theurl}\fi\par}

\vglue 7pt plus 0.3fil minus 3pt

{\bf Abstract}
\vglue 5pt plus 0.1fil minus 2pt

\theabstract

\vglue 7pt plus 0.3fil minus 3pt

{\bf AMS Classification numbers}\quad Primary:\quad \theprimaryclass

Secondary:\quad \thesecondaryclass

\vglue 5pt plus 0.3fil minus 2pt

{\bf Keywords}\quad \thekeywords

\vglue 10pt plus 0.5fil minus 5pt

{\small  Proposed: \theproposer\hfill Received: \receiveddate\nl
Seconded: \theseconders\hfill 
\ifx\reviseddate\relax                         
Accepted: \accepteddate                        
\else
Revised: \reviseddate                          
\fi}
\eject
}       

\let\maketitlepage\maketitlep
\let\maketitle\maketitlepage

\font\phead=cmsl9 scaled 950
\font=cmsl9 scaled 1050
\font\pnum=cmbx10 scaled 913
\font\lnum=cmbx10 
\font\pfoot=cmsl9 scaled 950
\font=cmsl9 scaled 1050
\ifplaintex
\headline{\vbox to 0pt{\vskip -4.5mm\line{\small\phead\ifnum
\count0=\startpage ISSN 1364-0380 (on line)
1465-3060 (printed) \hfill {\pnum\folio}\else\ifodd\count0\def\\{ }%
\ifx\theshorttitle\relax\thetitle\else\theshorttitle\fi\hfill{\pnum\folio}
\else\def\\{ and }{\pnum\folio}\hfill\ifx\theshortauthors\relax\theauthors
\else\theshortauthors\fi\fi\fi}\vss}}
\footline{\vbox to 0pt{\vglue 0mm\line{\small\pfoot\ifnum\count0=\startpage
\copyright\ \gtp\hfill\else
\gt, Volume \thevolumenumber\ (\thevolumeyear)\hfill\fi}\vss
}}
\else
\makeatletter
\def\@oddhead{{\small\lhead\ifnum\count0=\startpage ISSN 1364-0380 (on line)
1465-3060 (printed) \hfill {\lnum\number\count0}\else\ifodd\count0
\def\\{ }\ifx\theshorttitle\relax \thetitle \else\theshorttitle\fi\hfill
{\lnum\number\count0}\else\def\\{ and }{\lnum\number\count0}
\hfill\ifx\theshortauthors\relax 
\theauthors\else\theshortauthors\fi\fi\fi}}\def\@evenhead{\@oddhead}
\def\@oddfoot{\small\lfoot\ifnum\count0=\startpage\copyright\ \gtp\hfill\else
\gt, Volume \thevolumenumber\ (\thevolumeyear)\hfill\fi}
\def\@evenfoot{\@oddfoot}
\makeatother
\fi

\newwrite\gtoutfile
\long\gdef\makeheadfile{  
{\def\\{, }\def\s{ }
\immediate\openout\gtoutfile head.xxx
\immediate\write\gtoutfile{To: math@arxiv.org}
\immediate\write\gtoutfile{Subject: put or rep NNNNN:pppp}
\immediate\write\gtoutfile{--text follows this line--}
\immediate\write\gtoutfile{Proxy-for: \ifx\theasciiauthors\relax
\theauthors\else\theasciiauthors\fi\s<\ifx\theasciiemail\relax\theemail\else\theasciiemail\fi>}
\immediate\write\gtoutfile{\noexpand\\}
\immediate\write\gtoutfile{Authors: \ifx\theasciiauthors\relax
\theauthors\else\theasciiauthors\fi}
{\def\\{ }\immediate\write\gtoutfile{Title: \ifx\theasciititle\relax
\thetitle\else\theasciititle\fi}}
\immediate\write\gtoutfile{Subj-class: GT or SG or MG etc}
\immediate\write\gtoutfile{MSC-class: \theprimaryclass\ifx\thesecondaryclass\relax\else, \thesecondaryclass\fi}
\immediate\write\gtoutfile{Journal-ref: Geom. Topol. \thevolumenumber
(\thevolumeyear) \startpage-\finishpage}
\immediate\write\gtoutfile{Comments: Published by Geometry and Topology at}
\immediate\write\gtoutfile{\s\s http://www.maths.warwick.ac.uk/gt/GTVol\thevolumenumber/paper\thepapernumber.abs.html}
\immediate\write\gtoutfile{\noexpand\\}
\immediate\write\gtoutfile{}
\ifx\theasciiabstract\relax
\immediate\write\gtoutfile{\theabstract}\else
\immediate\write\gtoutfile{\theasciiabstract}\fi
\immediate\write\gtoutfile{}
\immediate\write\gtoutfile{\noexpand\\}
\immediate\write\gtoutfile{}
\immediate\closeout\gtoutfile}}  

\def\maketitlepage{\maketitlep\makeheadfile}
\let\maketitle\maketitlepage

\lognumber{225}

\volumenumber{6}
\papernumber{19}
\volumeyear{2002}
\pagenumbers{541}{562}
\received{15 January 2002}
\revised{9 August 2002}
\accepted{19 November 2002}
\published{24 November 2002}
\proposed{Walter Neumann}
\seconded{Cameron Gordon, Joan Birman}

\usepackage{amssymb,epsf}

\let\tilde\widetilde

\newtheorem{pro}{Proposition}[section]
\newtheorem{thm}[pro]{Theorem}
\newtheorem{lem}[pro]{Lemma}

\newtheorem{cnj}[pro]{Conjecture}
\newtheorem{cor}[pro]{Corollary}

\theoremstyle{definition}

\theoremstyle{remark}

\def\Ker{{\rm Ker}}
\def\Image{{\rm Image}}
\def\dim{{\rm dim}}
\def\genus{{\rm genus}}
\def\cycles{{\rm cycles}}

\begin{document}

\title{Virtual Betti numbers of genus 2 bundles}
\author{Joseph D Masters} 

\address{Mathematics Department, Rice University\\
Houston, TX 77005-1892, USA}

\email{mastersj@math.rice.edu}

\begin{abstract}
 We show that if $M$ is a surface bundle over $S^1$ with 
 fiber of genus 2,
 then for any integer $n$, $M$ has a finite cover $\tilde{M}$
 with $b_1(\tilde{M}) > n$.  A corollary is that $M$ can be
 geometrized using only the ``non-fiber" case of Thurston's
 Geometrization Theorem for Haken manifolds.
\end{abstract}

\asciiabstract{We show that if M is a surface bundle over S^1 with
 fiber of genus 2, then for any integer n, M has a finite cover
 tilde(M) with b_1(tilde(M)) > n.  A corollary is that M can be
 geometrized using only the `non-fiber' case of Thurston's
 Geometrization Theorem for Haken manifolds.}



\maketitlepage

\let\\\par

\section{Introduction}
Let $M$ be a 3--manifold.
Define the \textit{virtual first Betti number of $M$}
 by the formula
 $vb_1(M) = \sup \{ b_1(\tilde{M}) : \tilde{M}
 \textrm{ is a finite cover of M } \}$.

The following well-known conjecture is a strengthening of
 Waldhausen's conjecture about virtually Haken 3--manifolds.

\begin{cnj} \label{vb1conj}
Let $M$ be a closed irreducible 3--manifold with infinite fundamental group.
  Then either $\pi_1 M$ is virtually solvable, or $vb_1(M) = \infty$.
\end{cnj}

 Combining the Seifert Fiber Space Theorem, the Torus Theorem, and
 arguments involving characteristic submanifolds, 
 Conjecture \ref{vb1conj} is known to be true
 in the case that $\pi_1 M$
 contains a subgroup isomorphic to $\mathbb{Z} \times \mathbb{Z}$.
 However, little is known in the atoroidal case.

 In \cite{G}, Gabai called attention to Conjecture \ref{vb1conj} in the
 case that $M$ fibers over $S^1$. This seems a natural place to start, 
 in light of Thurston's conjecture that every closed hyperbolic 3--manifold
 is finitely covered by a bundle.   The purpose of this paper is to give
 some affirmative results for this case.  In particular, we prove
 Conjecture \ref{vb1conj} in the case where $M$ is a genus 2 bundle.

 Throughout this paper, if $f\co  F \rightarrow F$ is an automorphism of
 a surface, then $M_f$ denotes the associated mapping torus.
 Our main theorem is the following:

\begin{thm} \label{main}
 Let $f \co F \rightarrow F$ be an automorphism of a surface.
  Suppose there is a finite group $G$
 of automorphisms of $F$, so that $f$ commutes with each element of $G$,
 and $F / G$ is a torus with at least one cone point.
Then $vb_1(M_f) = \infty$.
\end{thm}

We have the following corollaries:

\begin{cor} \label{hyperelliptic}
Suppose $F$ has genus at least 2, and $f\co F \rightarrow F$
 is an automorphism which commutes with a hyper-elliptic
 involution on $F$. Then $vb_1(M_f) = \infty$.
\end{cor}

\begin{proof}
 Let $\tau$ be the hyper-elliptic involution.
  Since $f$ commutes with $\tau$, $f$ induces an automorphism 
 $\bar{f}$ of $F/ \tau$,
 which is a sphere with $2g+2$ order 2 cone points.
 $F/ \tau$ is double covered by a hyperbolic orbifold $T$,
 whose underlying space is a torus.  
 By passing to cyclic covers of $M$, we may replace $f$ (and $\bar{f}$)
 with powers, and so we may assume $\bar{f}$ lifts to $T$.
 Corresponding to $T$, there is a 2--fold cover $\tilde{F}$ of $F$
 to which $f$ lifts, and an associated cover $\tilde{M}$ of $M$
 whose monodromy satisfies the hypotheses of Theorem \ref{main}. 
 \end{proof}

\begin{cor} \label{genus2}
Let $M$ be a surface bundle with fiber $F$ of genus 2.
 Then $vb_1(M) = \infty$.
\end{cor}

\begin{proof}
 Since the fiber has genus 2, the monodromy map commutes 
(up to isotopy) with the central hyper-elliptic involution on $F$.
The result now follows from Corollary \ref{hyperelliptic}.
\end{proof}

 To state our next theorem, we require some notation.
 Recall that, by \cite{H}, the mapping class group of a surface is generated
 by Dehn twists in the loops pictured in Figure 1.
 If $\ell$ is a loop in a surface, we let $D_{\ell}$
 denote the right-handed Dehn twist in $\ell$.

\begin{figure}[ht!]
\cl{\epsfxsize.95\hsize\epsfbox{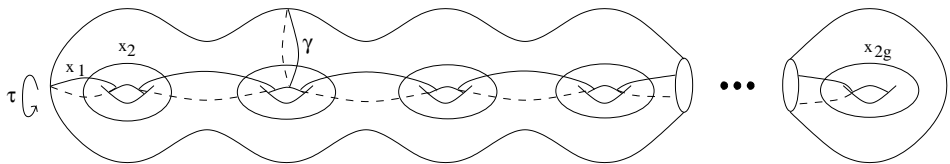}}
\caption{The mapping class group is generated by Dehn twists in these
 loops.}
\end{figure}

 With the exception of $D_{\gamma}$, these Dehn twists
 each commute with the involution $\tau$ pictured in Figure 1.
 Let $H$ be the subgroup of the mapping class group
 generated by the $D_{x_i}$'s.
 For any monodromy $f \in H$, we may apply Corollary \ref{hyperelliptic}
 to show that the associated bundle $M$ has $vb_1(M) = \infty$.
 The proof provides an explicit construction of covers--
 a construction which may be applied to any bundle, regardless of monodromy. 
 These covers will often have extra
 homology, even when the monodromy does not commute with $\tau$.
 For example, we have the following theorem, which is proved in Section 7.

\begin{thm} \label{nonhyper}
Let $M$ be a surface bundle over $S^1$ with fiber $F$
 and monodromy $f\co  F \rightarrow F$.
 Suppose that $f$ lies in the subgroup of the mapping class
 group generated by $D_{x_1}, ..., D_{x_{2g}}$ and $D_{\gamma}^8$.
 Then $vb_1(M) = \infty$.
\end{thm}

None of the proofs makes any use of a geometric structure.  In
 fact, for a bundle satisfying the hypotheses of one
 of the above theorems, we may give
 an alternative proof
 of Thurston's hyperbolization theorem for fibered 3--manifolds.
For example, we have:
\begin{thm}[Thurston]
Let $M$ be an atoroidal surface bundle over $S^1$
 with fiber a closed
 surface of genus 2.  Then $M$ is hyperbolic.
\end{thm} 
\begin{proof}
By Corollary \ref{genus2}, $M$ has a finite cover $\tilde{M}$
 with $b_1(\tilde{M}) \geq 2$.
 Therefore, by \cite{T},
 $\tilde{M}$
 contains a non-separating incompressible surface which is not
 a fiber in a fibration.
 Now the techniques of the non-fiber case of Thurston's
 Geometrization Theorem (see \cite{O}) may be applied
 to show that $\tilde{M}$ is hyperbolic.
 Since $M$ has a finite cover which is hyperbolic,
 the Mostow Rigidity Theorem implies that $M$
 is homotopy equivalent to a hyperbolic 3--manifold.
 Since $M$ is Haken, Waldhausen's Theorem (\cite{W})
 implies that $M$ is in fact homeomorphic
 to a hyperbolic 3--manifold.
\end{proof}
 
 We say that a surface automorphism $f \co  F \rightarrow F$
 is \textit{hyper-elliptic} if it commutes with some hyperelliptic
 involution on $F$.  Corollary \ref{hyperelliptic}
 prompts the question:  is a hyper-elliptic monodromy
 always attainable in a finite cover?  Our final theorem
 shows that the answer is no.

\begin{thm} \label{badex}
  There exists a closed surface $F$, and a
  pseudo-Anosov automorphism $f\co  F \rightarrow F$, such that
 $f$ does not lift to become hyper-elliptic in any finite
 cover of $F$.
\end{thm}

The proof of Theorem \ref{badex} will be given in Section 8.

\rk{Acknowledgements}
 I would like to thank Andrew Brunner, Walter Neumann
 and Hyam Rubinstein, whose work suggested the
 relevance of punctured tori to this problem.
 I also thank Mark Baker, Darren Long, Alan Reid and the referee
 for carefully reading previous versions
 of this paper, and providing many helpful comments.
 Alan Reid also helped with the proof of Theorem \ref{badex}.
 Thanks also to The University
 of California at Santa Barbara, where this work was begun.

This research was supported by an NSF Postdoctoral Fellowship.

\section{Homology of bundles: generalities}

 In what follows, we shall try to keep notation to a minimum;
 in particular we shall often neglect to distinguish notationally
 between the monodromy map $f$, and the various maps
 which $f$ induces on covering spaces or projections.
 All homology groups will be taken with $\mathbb{Q}$ coefficents.

 Suppose $f$ is an automorphism of a closed 2--orbifold $O$.
 The mapping torus $M_f$ associated with $O$ is a 3--orbifold,
 whose singular set is a link.  
 We have the following well-known formula for the first Betti number
 of $M_f$:
\begin{equation} \label{b1}
b_1(M_f) = 1 + \dim({\rm fix}(f_*)),
\end{equation}
 where ${\rm fix}(f_*)$ is the subspace of $H_1(O)$ on which $f_*$ acts trivially.
 This can be derived by abelianizing the standard HNN presentation
 for $\pi_1 M_f$.

 Suppose now that $O$ is obtained from a punctured surface $F$ by filling in
 the punctures with disks or cone points, and suppose that
 $f$ restricts to an automorphism of $F$.  
 Let $V \subset H_1(F, \partial F)$
 be the subspace on which the induced map $f_*$ acts trivially.
 Then the first Betti number for the mapping torus of $O$
 can also be computed by the following formula.

\begin{pro} \label{b1orb}
For $M_f$ and $V$ as above, we have  $b_1(M_f) = 1 + \dim(V)$.
\end{pro}

\begin{proof}
By Formula \ref{b1}, $b_1(M_f) = 1 + \dim(W)$,
 where $W \subset H_1(O)$ is the subspace
 on which $f_*$ acts trivially.
 
 Let $i\co  F \rightarrow O$ be the inclusion map,
 and let $K$ be the kernel
 of the quotient map from $H_1(F)$ onto $H_1(F, \partial F)$.
 The cone-point relations imply that every element in
 $i_* K$ is a torsion element in $H_1(O)$; 
 since we are using $\mathbb{Q}$--coefficients,
 $i_* K$ is in fact trivial in $H_1(O)$.
 The action of $f_*$ on $H_1(O)$ is therefore identical to the
 action of $f_*$ on $H_1(F, \partial F)$,
 so $\dim(W) = \dim(V)$, which proves the formula. 
\end{proof}

We will also need the following technical proposition.

\begin{pro} \label{pullback}
Let $F$ be a punctured surface, and let $f\co F \rightarrow F$
 be an automorphism which fixes the punctures.
  Let $F^+$ be a surface obtained
 from $F$ by filling in one or more of the punctures,
 and let $f^+\co F^+ \rightarrow F^+$ be the map induced by $f$.
 Suppose $\widetilde{F^+}$ is a cover of $F^+$,
 such that $f^+$ lifts, and suppose
 $\alpha^+ \in \widetilde{F^+}$ is a loop which misses all filled-in punctures,
 and such that
 $f^+ [\alpha^+]=[\alpha^+]\in H_1(\widetilde{F^+},\partial \widetilde{F^+})$.
 Let $\widetilde{F}$ be the cover of $F$ corresponding
 to the cover $\widetilde{F^+}$ of $F^+$, and let
 $i\co  \widetilde{F} \rightarrow \widetilde{F^+}$ be the natural
 inclusion map.  Let $\alpha = i^{-1} \alpha^+$.
 Then $f [\alpha] = [\alpha] \in H_1(\widetilde{F}, \partial \widetilde{F})$.
\end{pro}

\begin{proof}
 The surface $\widetilde{F^+}$ is obtained from $\widetilde{F}$ by
 filling in a certain number of punctures, say
 $\beta_1, ..., \beta_k$, of $\widetilde{F}$.
 The map
$f\co  H_1(\tilde{F}, \partial \tilde{F}) \rightarrow H_1(\tilde{F}, \partial \tilde{F})$
 may be obtained from the map $f\co  \pi_1 \tilde{F} \rightarrow \pi_1 \tilde{F}$
 by:
\medskip

(1)\qua First add the relations $\beta_1, ..., \beta_k = id$.
There is an induced map
 $$f\co  \pi_1 \tilde{F}/<\beta_1= ... = \beta_k = id>
 \rightarrow \pi_1 \tilde{F}/<\beta_1= ... = \beta_k = id>.$$

(2)\qua Add the relations which kill the remaining boundary components.
\medskip

(3)\qua Add the relations $[x,y] = id$ for all
 $x, y \in \pi_1 \widetilde{F^+}$.\\
\medskip

After completing step 1, one has precisely the action
 of $f$ on $\pi_1 \widetilde{F^+}$.  After completing steps 2
 and 3, one then has the action
 of $f$ on $H_1(\widetilde{F^+}, \partial \widetilde{F^+})$.
 So the action of $f$ on these groups is identical,
 and $[\alpha]$ is a fixed class.
\end{proof}

If $\Gamma$ is a group, we may define $b_1(\Gamma)$ to be the
 $\mathbb{Q}$--rank of its abelianization, and the virtual first Betti number
 of $\Gamma$ by
$$vb_1(\Gamma) = \sup \{ b_1(\tilde{\Gamma}) : \tilde{\Gamma} \textrm{ is a finite index subgroup of } \Gamma \}.$$
Clearly, for a 3--manifold $M$, $vb_1(M) = vb_1(\pi_1(M))$.
We have the following:

\begin{lem} \label{vb1homo}
Suppose $\Gamma$ maps onto a group $\Delta$.  Then
 $$(vb_1(\Gamma) - b_1(\Gamma)) \geq (vb_1(\Delta) - b_1(\Delta)).$$
\end{lem}

Before proving this, we will need a preliminary lemma.
 We let $H_1(\Gamma)$ denote the abelianization of $\Gamma$,
 tensored over $\mathbb{Q}$.
 Representing $\Gamma$ by a 2--complex $C_{\Gamma}$, then
 $H_1(\Gamma) \cong H_1(C_{\Gamma})$.

 Any subgroup $\widetilde{\Gamma}$ of $\Gamma$ determines a 2--complex
 $\widetilde{C_{\Gamma}}$ and a covering map
 $p \co \widetilde{C_{\Gamma}} \rightarrow C_{\Gamma}$.
 We can define
 a map $j\co H_1(C_{\Gamma}) \rightarrow H_1(\widetilde{C_{\Gamma}})$
 by the rule $j([\ell]) = [p^{-1} \ell]$,
 for any loop $\ell$ in $C_{\Gamma}$.
 If $\ell$ bounds a 2--chain in $C_{\Gamma}$,
 then $p^{-1} \ell$ bounds a 2--chain in $\widetilde{C_{\Gamma}}$,
 so this map is well-defined.  Using the isomorphisms between the homology
 of the groups and the homology of the 2--complexes, we get a map,
 which we also call $j$, from $H_1(\Gamma)$ to $H_1(\tilde{\Gamma})$.

\begin{lem} \label{inject}
If $\tilde{\Gamma}$ has finite index, then the map $j$ is injective.
\end{lem}

\begin{proof}
Suppose $[\gamma] \in \Ker(j)$,
 where $\gamma$ is an element of $\Gamma$, and let $\ell \in C_{\Gamma}$
 be a corresponding loop.
 Then $[\ell] \in \Ker(j)$, so $[p^{-1} \ell] = 0$, and therefore
 $$0 = p_*[p^{-1} \ell] = n[\ell],$$
 where $n$ is the index of $\widetilde{\Gamma}$. 
 Since we are using $\mathbb{Q}$--coefficients,
 $H_1(C_{\Gamma})$ is torsion-free, so $[\ell] = 0$ in $H_1(C_{\Gamma})$, and
 therefore $[\gamma] = 0$ in $H_1(\Gamma)$.
\end{proof}

\begin{proof}[Proof of Lemma \ref{vb1homo}]
Let $f\co  \Gamma \rightarrow \Delta$ be a surjective map.
We have the following commutative diagram:
$$\begin{array}{lllcl}
\; & \;          & \;             & i_1         & \; \\
1  & \longrightarrow & \widetilde{\Gamma} & \longrightarrow & \Gamma   \\
\; &    \;       & \downarrow g    & i_2         & \downarrow f\\
1  & \longrightarrow & \widetilde{\Delta} & \longrightarrow & \Delta,     
\end{array}$$
where $i_1$ and $i_2$ are inclusion maps,
 and the surjective map $g$ is induced from the other maps.
There is an induced diagram on the homology:
$$\begin{array}{ccc}
 \; & i_{1*} & \; \\
H_1(\widetilde{\Gamma})
      & \longrightarrow & H_1(\Gamma) \\
 \downarrow g_*    & i_{2*}          & \downarrow f_*\\
H_1(\tilde{\Delta}) & \longrightarrow
 & H_1(\Delta).
\end{array}$$
Let $j_1\co  H_1(\Gamma) \longrightarrow H_1(\tilde{\Gamma})$
 and $j_2 \co  H_1(\Delta) \longrightarrow H_1(\tilde{\Delta})$ be
 the injective maps given by Lemma \ref{inject}. 
 These maps give rise to the following diagram, which
 can be checked to be commutative:
$$\begin{array}{cccll}
 \; & j_1 & \; \\
H_1(\tilde{\Gamma})      & \longleftarrow & H_1(\Gamma) & \longleftarrow & 0 \\
 \downarrow g_*    & j_2          & \downarrow f_*\\
H_1(\tilde{\Delta}) & \longleftarrow & H_1(\Delta) & \longleftarrow & 0.
\end{array}$$
The definitions of the maps give
 that
$$i_{1*}j_1( [\gamma]) = n[\gamma],\leqno{(*)}$$
 and a similar relation
 for $i_{2*}$ and $j_2$.
Therefore $\Ker(i_{1*})$ and $\Image(j_1)$
 are disjoint subspaces of $H_1(\tilde{\Gamma})$. 
 Also,
\begin{eqnarray*} 
\dim(H_1(\tilde{\Gamma})) &=& \dim(\Ker(i_{1*})) + \dim(\Image(i_{1*}))\\
 &=& \dim(\Ker(i_{1*})) + \dim(H_1(\Gamma)), \textrm{ by the relation (*)}\\
 &=& \dim(\Ker(i_{1*})) + \dim(\Image(j_1)),\\
\end{eqnarray*}
 so we get $H_1(\tilde{\Gamma}) = \Ker(i_{1*}) \times \Image(j_1)$,
 and similarly $H_1(\tilde{\Delta}) = \Ker(i_{2*}) \times \Image(j_2)$.
 Substituting these decompositions into the previous diagram gives:
$$\begin{array}{cccll}
 \; & j_1 & \; \\
\Ker(i_{1*}) \times \Image(j_1) 
      & \longleftarrow & H_1(\Gamma) & \longleftarrow & 0 \\
 \downarrow g_*    & j_2          & \downarrow f_*\\
\Ker(i_{2*}) \times \Image(j_2) & \longleftarrow & H_1(\Delta) & \longleftarrow & 0.
\end{array}$$
 By the commutativity of this diagram, we have that
 $g_* \Image(j_1) \subset \Image(j_2)$.
 Also, by the commutativity of a previous
 diagram, we have $g_* \Ker(i_{1*}) \subset \Ker(i_{2*})$.
 Since $g_*$ is surjective, we must therefore have
 $g_* \Ker(i_{1*}) = \Ker(i_{2*})$, so
\newline
 $\dim(\Ker(i_{1*})) \geq \dim(\Ker(i_{2*}))$,
 from which the lemma follows.
\end{proof}

\section{Reduction to a once-punctured torus}

We are given an automorphism of a torus with an arbitrary number, $k$,
 of cone points.  We denote this orbifold $T(n_1, ..., n_k)$,
 where $n_i$ is the order of the $i$-th cone point.
 Let $\mathcal{M}(T(n_1, ..., n_k))$ be the mapping class group of
 $T(n_1, ..., n_k)$.  In general, these groups are rather complicated.
 However, the mapping class
 group of a torus with a single cone point is quite simple, being
 isomorphic to $SL_2(\mathbb{Z})$.

 Let $\mathcal{M}_0(T(n_1, ..., n_k))$ denote the finite-index
 subgroup of the mapping class group which consists of those automorphisms
 which fix all the cone points of $T(n_1, ..., n_k)$.
 The following elementary fact allows us to pass to the simpler
 case of a single cone point. 

\begin{lem}
For any $i$, there is a homomorphism
 $\theta_i\co \mathcal{M}_0(T(n_1, ..., n_k))$ onto $\mathcal{M}(T(n_i))$.
\end{lem}

\begin{proof}
Let $f \in \mathcal{M}_0(T(n_1, ..., n_k))$.  Since $f$ fixes
 the cone points, it restricts to a map on the punctured surface
 which is the complement of all the cone points except the $i$th one.
 After filling in these punctures, there is an induced map
 $\theta_i(f)$ on $T_{n_i}$.  It is easy to see that this is well-defined,
 surjective, and a homomorphism.
\end{proof}

\begin{lem} \label{bundlehomo}
Let $f \in \mathcal{M}_0(T(n_1, ..., n_k))$.
 Then there is a surjective homomorphism from
 $\pi_1 M_f \rightarrow \pi_1 M_{\theta_i f}$.
\end{lem}

\begin{proof}
Let $F$ be the punctured surface obtained from $T(n_1, ..., n_k)$
 by removing all the cone points.
 Let $x_1, ..., x_k \in \pi_1 F$ be loops around
 the $k$ cone points, and complete
 these to a generating set with loops $x_{k+1}, x_{k+2}$.
 We have:
\begin{eqnarray*}
\pi_1 M_f \cong <x_1,..., x_{k+2}, t>/<t x_1 t^{-1} &=& f x_1,
                \,\, ... \,\, , \,\, t x_{k+2} t^{-1} = f x_{k+2},\\
             x_1^{n_1} &=& ... \,\, =x_k^{n_k} = 1>.
\end{eqnarray*}
 From this presentation a presentation for $ \pi_1 M_{\theta_i f}$
 may be obtained by adding the additional relations
 $x_j = id$, for all $j \leq k, j \neq i$.   
\end{proof}

\begin{cor} \label{fibundlehomo}
Let $f \in \mathcal{M}(T(n_1, ..., n_k))$.  Then there
 is a finite index subgroup of $\pi_1 M_f$
 which maps onto $\pi_1 M_g$, where $g$ is an
 automorphism of a torus with a single cone point.
\end{cor}
 
\begin{proof}
By passing to a finite-index subgroup, we may replace $f$ with
 a power, and then apply Lemma \ref{bundlehomo}.
\end{proof}

\section{Increasing the first Betti number by at least one}

Before proving Theorem \ref{main}, we first prove:

\begin{lem} \label{vb1>1}
Let $M_f$ be as in the statement of Theorem \ref{main}. Then\\
$vb_1(M_f) > b_1(M_f)$.
\end{lem}

 We remark that this result, combined with Lemma \ref{vb1homo}
 and the arguments in the proof of Cor \ref{genus2}, implies
 that the first Betti number of a genus 2 bundle can be increased by
 at least 1.  

By Corollary \ref{fibundlehomo},
 Lemma \ref{vb1>1} will follow from the following lemma.

\begin{lem} \label{onepuncture}
Let $f \in \mathcal{M}(T(n))$ be an automorphism
 of a torus with a single cone point.  Then $M_f$ has
 a finite cover $\widetilde{M_f}$ such that
 $b_1(\widetilde{M_f}) > b_1(M_f)$.
\end{lem}

\begin{proof}[Proof of Lemma \ref{onepuncture}]

We shall use $T$ to denote the once-punctured torus
 obtained by removing the cone point
 of $T(n)$. There is an induced map $f \co  T \rightarrow T$.
 In order to construct covers of $T$, we require the
 techniques of \cite{M}.  For convenience, the relevant ideas
 and notations are contained in the appendix.
 In what follows, we assume familiarity with this material.

\medskip

{\bf Case 1}\qua $n=2$

We let $J$ denote the subgroup of the mapping
 class group of $T$ generated by $D_x$ and $D_y^4$.
 By Lemma \ref{fi}, $J$ has finite index, so
 we may assume, after replacing $f$ with a power,
 that the map  $f\co T \rightarrow T$ lies in $J$.
 
 As explained in the appendix, any four permutations
 $\sigma_1, ..., \sigma_4$ on $r$ letters will determine
 a $4r$--fold cover $\tilde{T}$ of $T$.
 We set:
\\
{\bf I}\quad $\sigma_2 = \sigma_1^{-1}$ and $\sigma_4 = \sigma_3^{-1}$,\\
\\
so $f$ lifts to $\tilde{T}$ by Lemma \ref{int0}.
We shall require every lift of $\partial T$ to unwrap once or twice
 in $\tilde{T}$. 
This property is equivalent to the following:\\
\\
{\bf II}\quad $(\sigma_i \sigma_{i+1}^{-1})^2 = id$ for all $i$.\\
\\
 To find permutations satisfying I and II, we consider the abstract group
 generated by the symbols $\sigma_1, ..., \sigma_4$, satisfying
 relations I and II.
 If this group surjects a finite group $G$, then we may take the
 associated permutation representation, and obtain
 permutations $\sigma_1, ..., \sigma_4$ on $|G|$ letters
 satisfying I and II.
 In the case under consideration, we may take $G$ to be a cyclic
 group of order 4.
 This leads to the representation $\sigma_1 = \sigma_3 = (1234)$,
 $\sigma_2 = \sigma_4 = \sigma_1^{-1}$.  The associated
 cover is pictured in Figure 2.
 
\begin{figure}[ht!]
\cl{\epsfxsize 3.5in \epsfbox{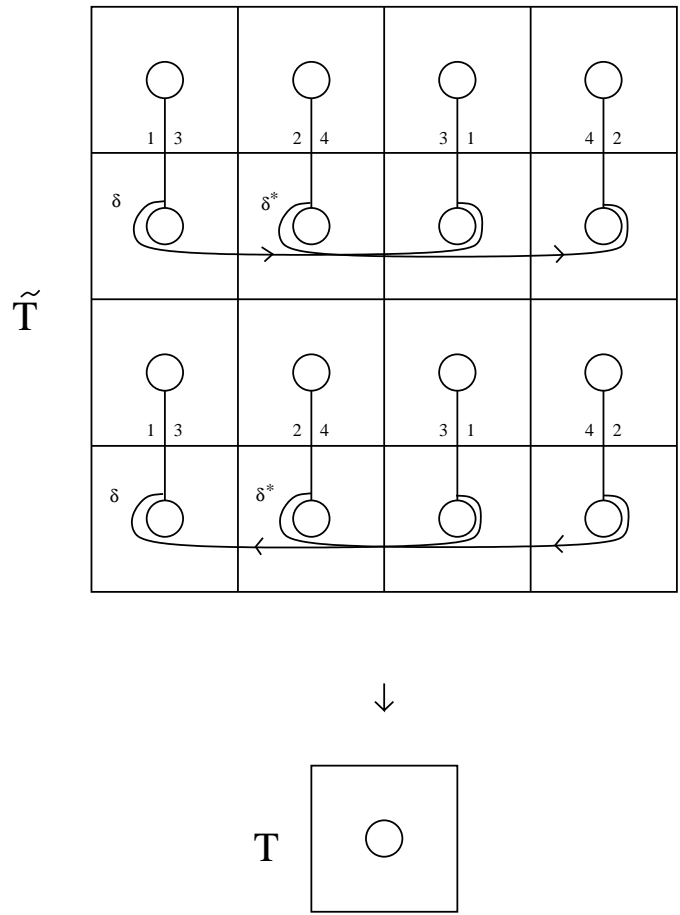}}
\caption{The cover $\tilde{T}$ of $T$}
\end{figure}

 Lemma \ref{int0} now guarantees that non-trivial
 fixed classes in $H_1(\tilde{T}, \partial \tilde{T})$ exist.
 Rather than invoke the lemma, however, we shall give the explicit
 construction for this simple case.
 Consider the classes
 $[\delta], [\delta^*] \in H_1(\tilde{T}, \partial \tilde{T})$
 pictured in Figure 2.

\begin{pro} \label{delta}
 $[\delta], [\delta^*] \in H_1(\tilde{T}, \partial \tilde{T})$
 are non-zero classes which are fixed by any element of $J$.
\end{pro}

In the proof, the notation $I(.,.)$ stands for the algebraic
 intersection pairing on $H_1$ of a surface.

\begin{proof}
 The fact that $[\delta]$ and $[\delta^*]$
 are non-peripheral follows from the fact that
 $I([\delta], [\delta^*]) = 2$.
 The loops $\delta$ and $\delta^*$
 have algebraic intersection number 0 with each lift of $y$,
 and therefore their homology classes are fixed by the lift of $D_y^4$.
 
 By Property I and by Lemma \ref{lift},
 $D_x$ lifts to $\tilde{T}$,
 and acts as the identity on Rows 2 and 4.
 In particular, $[\delta]$ and $[\delta^*]$ are
 fixed by the lift of $D_x$.
 Therefore,  $[\delta]$ and $[\delta^*]$ are fixed by
 any element of $J$.
\end{proof}

 Since $\partial T$ unwraps exactly twice in every
 lift to $\tilde{T}$, then by filling in the punctures of
 $\tilde{T}$, we obtain a manifold cover $\widetilde{T(n)}$
 of $T(n)$.
 Since $f$ lifts to $\tilde{T}$, then $f$ lifts to $\widetilde{T(n)}$.
 An application of Propositions \ref{delta} and \ref{b1orb}
 finishes the proof of Lemma \ref{onepuncture} in this case.

\medskip

{\bf Case 2}\qua $n \geq 3$
\\

In this case, we shall require a cover of $T$ in which
the boundary components unwrap $n$ times.
 We construct a cover $\tilde{T}$ of $T$,
 mimicking the construction given in Case 1.
 We start with the standard
 $\mathbb{Z}/r \times \mathbb{Z}/4$ cover of $T$,
 and alter it by cutting and pasting in a manner specified
 by permutations $\sigma_1, ..., \sigma_4$.
 By raising $f$ to a power, we may assume that
 $f$ lies in $J$, the subgroup of the mapping class group
 of $T$ generated by $D_x$ and $D_y^4$.
 Again, $D_y^4$ lifts to Dehn twists in the lifts of $y$.
 Lemma \ref{lift} shows that, if we set: \\
\\
{\bf I}\quad $\sigma_2 = \sigma_1^{-1}$ and $\sigma_4 = \sigma_3^{-1}$,\\
\\
 then $D_x$ lifts also, so $f$ lifts.
To ensure that the boundary components unwrap appropriately,
 we also require $(\sigma_i \sigma_{i+1}^{-1})^n = 1$.
 Combining this with condition I gives:\\
\\
{\bf II}\quad $(\sigma_1)^{2n} = (\sigma_3)^{2n} = (\sigma_1\sigma_3)^n = 1$.\\
\\
If we consider the symbols $\sigma_1$ and $\sigma_3$ as representing
 abstract group elements, Conditions I and II determine a hyperbolic triangle
 group $\Gamma$. It is a well-known property of triangle
 groups that $\sigma_1, \sigma_3$, and $\sigma_1\sigma_3$ will have
 orders in $\Gamma$ as given by the relators in Condition I.
 Also, it is a standard fact that in this case $\Gamma$ is infinite, and
 residually finite.  Therefore, 
 $\Gamma$ surjects arbitrarily large finite groups such that the images of
 $\sigma_1, \sigma_3$, and $\sigma_1\sigma_3$ have orders
 $2n, 2n$, and $n$, respectively. Let $G$ be such a finite
 quotient, of order $N$ for some large number $N$.
 By taking the left regular permutation representation of $G$,
 we obtain permutations on $N$ letters satisfying Conditions I
 and II, as required.
 
Let $V$  denote the subspace of
 $H_1(\tilde{T}, \partial \tilde{T})$ fixed by $f$. 
 By Lemma \ref{int0}, 
 $\dim(V) \geq 2\genus(R_2)$,
 where $R_2$ is the subsurface of $\tilde{T}$ corresponding to Row 2.

The formula for genus is:
$$\genus(R_2) = \hbox{$1\over2$}(2 - \chi(R_2) - (\# \textrm{ of punctures of }R_2)).$$
 Any permutation $\sigma$ decomposes uniquely as a product
 of disjoint cycles; we denote the set of these cycles by
 $\cycles(\sigma)$.
 The punctures of $R_2$
 are in 1--1 correspondence with the cycles of $\sigma_1$, $\sigma_3$ and
 $\sigma_3\sigma_1$.
 Also, since $R_2$ is an $N$--fold cover of a thrice-punctured sphere, we
 have the Euler characteristic $\chi(R_2) = -N$,
 and so we get:
 $$ \genus(R_2) = \hbox{$1\over2$}( 2+ N - (|\cycles(\sigma_1)| + |\cycles(\sigma_3)|
 + |\cycles(\sigma_1\sigma_3)|)).$$
 Recall that an \textit{m--cycle} is a permutation which is conjugate
 to $(1 ... m)$.
 Any permutation $\sigma$ coming from the left regular permutation
 representation of $G$ decomposes as a product of
 $N/order(\sigma)$ disjoint $order(\sigma)$--cycles,
 and therefore
$$|\cycles(\sigma_1 \sigma_3)| =
  2|\cycles(\sigma_1)| = 2|\cycles(\sigma_3)| = |G|/n = N/n.$$
Combining the above formulas gives
$$\dim(V) \geq  2 + N(1 - 2/n).$$
 So $\dim(V)$ can be made arbitrarily large.

 There are corresponding covers $\tilde{T}$ of $T$,
 and $\widetilde{T(n)}$ of $T(n)$.
 Proposition \ref{b1orb} then
 shows that, in this case, $vb_1(M_f) = \infty$.
\end{proof}

\section{Infinite virtual first Betti number}

In this section we prove Theorem \ref{main}.

\begin{lem} \label{vb1=infty}
Let $f \in \mathcal{M}(T(n))$ be an automorphism
 of a torus with a single cone point.  Then $vb_1(M_f) = \infty$.
\end{lem}

\begin{proof}

 In the course of proving Lemma \ref{vb1>1}, we actually proved
 Lemma \ref{vb1=infty} in the case $n > 2$,
 so we assume $n = 2$.  The proof of Lemma \ref{vb1>1} also shows
 how to increase $b_1(M_f)$ by at least 2;  next we will show how
 to increase $b_1(M_f)$ by at least 4, and fianlly we will indicate
 how to iterate this process to increase $b_1(M_f)$ arbitrarily.
 
 Again, let $T$ be the once punctured torus obtained by removing
 the cone point of $T(2)$.
 By replacing $T$ with a 2--fold cover,
 and by replacing $f$ with a power (to make it lift),
 we may assume that $T$ has two boundary components,
 denoted $\alpha_1, \alpha_2$.
 By again replacing $f$ with a power, we may assume
 that $f$ fixes both $\alpha_i$'s.

 Let $T_1^+$ denote the once punctured torus obtained by filling in
 $\alpha_2$.
 Since $f$ fixes both $\alpha_i$'s, there is an induced
 automorphism $f\co  T_1^+ \rightarrow T_1^+$.
 Let $\widetilde{T_1^+}$ be the 16--fold cover of $T_1^+$
 as constructed in the previous section,
 and let $\delta_1^+, \delta_1^{+*}$ be the loops constructed
 previously, whose homology classes are fixed by (a power of) $f$.

 Let $\widetilde{T}_1$ denote the cover of $T$ corresponding to
 $\widetilde{T_1^+}$ (see Figure 3).
 By replacing $f$ with a power, we may assume
 that $f$ lifts to $\widetilde{T_1}$.
 Let $\delta_1, \delta_1^* \subset \widetilde{T}_1$
 denote the pre-images of $\delta_1^+$ and $\delta_1^{+*}$ under
 the natural inclusion map
 (after an isotopy, we may assume that
 $\delta_1^+$ and $\delta_1^{+*}$ are disjoint
 from all filled-in punctures, so that
 $\delta_1$ and $\delta_1^*$ are in fact loops).

 Since $[\delta_1^+]$ and $[\delta_1^{+*}]$ are fixed
 classes in $H_1(\widetilde{T}_1^+, \partial \widetilde{T}_1^+)$, then by
 Proposition \ref{pullback},
 $[\delta_1]$ and $[\delta_1^*]$ are fixed classes
 in $H_1(\widetilde{T}_1, \partial \widetilde{T}_1)$.
 Note that $I([\delta_1], [\delta_1^*]) = 2$.

 Starting with $\alpha_2$ instead of $\alpha_1$, we may perform 
 the analogous construction to obtain a cover
 $\widetilde{T_2}$ of $T$ containing fixed classes
 $[\delta_2], [\delta_2^*] \in H_1(\widetilde{T}_2, \partial \widetilde{T}_2)$,
 with algebraic intersection number 2.
 Moreover, as indicated by Figure 3,
 $\delta_2$ and $\delta_2^*$ may be chosen so that their
 projections to $T$ are disjoint from the projections of
 $\delta_1$ and $\delta_1^*$ to $T$.
 
\begin{figure}[ht!]
\cl{\epsfxsize\hsize\epsfbox{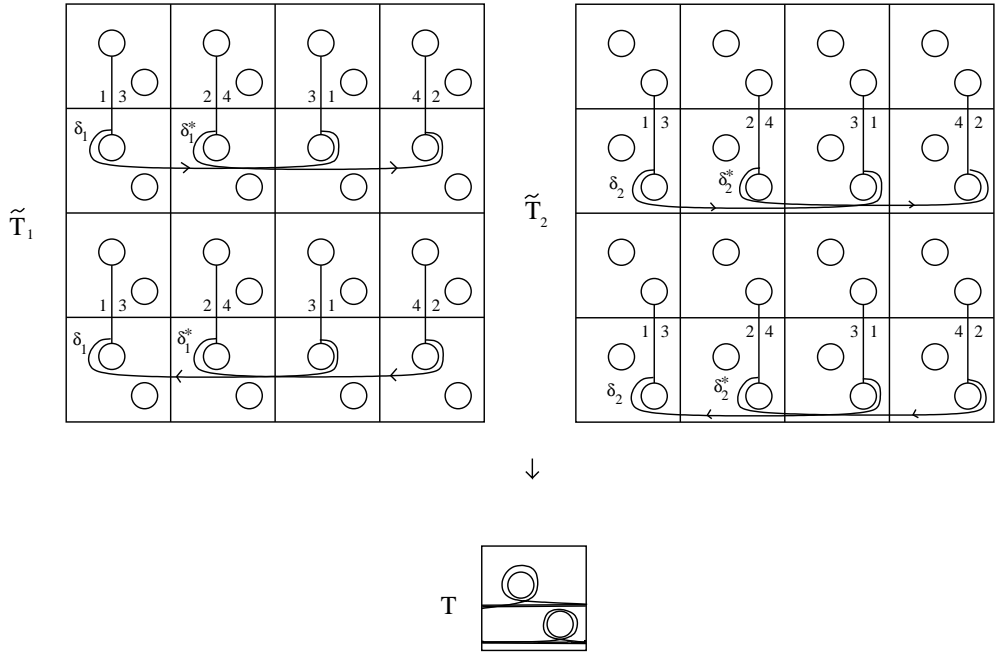}}
\caption{We may arrange for $\delta_1 \cup \delta_1^*$
 and $\delta_2 \cup \delta_2^*$ to have disjoint projections.}
\label{cover2}
\end{figure}

 Let $\widetilde{T}$ denote the cover of $T$
 with covering group
 $\pi_1(\widetilde{T_1}) \cap \pi_1(\widetilde{T_2})$. 
 Since $f$ lifts to $\widetilde{T_1}$
 and $\widetilde{T_2}$, then $f$ also lifts to $\widetilde{T}$.
 Let $\tilde{\delta}_i$ and $\tilde{\delta}_i^*$ denote the
 full pre-images in $\widetilde{T}$ of $\delta_i$
 and $\delta_i^*$, respectively.
 Recall that by construction,
 $\delta_i$ and $\delta_i^*$ have the following properties,
 for $i = 1, 2$:

\medskip
(1)\qua $[\delta_i], [\delta_i^*] \in H_1(\widetilde{T}_i, \partial \widetilde{T}_i)$ 
 are fixed classes.\\\medskip
(2)\qua $I([\delta_i], [\delta_i^*]) \neq 0$.\\\medskip
(3)\qua The projections of $\delta_1 \cup \delta_1^*$ and $\delta_2 \cup \delta_2^*$
 to $T$ are disjoint.\\\medskip 
\\
 Therefore, by elementary covering space arguments, we deduce that
 $\tilde{\delta}_i$ and $\tilde{\delta}_i^*$ have the following
 properties for $i = 1, 2$:\\\medskip
\\
(1)\qua $[\tilde{\delta}_i], [\tilde{\delta}_i^*]
 \in H_1(\widetilde{T}, \partial \widetilde{T})$ are fixed  classes.\\\medskip
(2)\qua $I([\tilde{\delta}_i], [\tilde{\delta}_i^*]) \neq 0$.\\\medskip
(3)\qua $\tilde{\delta}_1 \cup \tilde{\delta}_1^*$ and
 $\tilde{\delta}_2 \cup \tilde{\delta}_2^*$ are disjoint.\\

\medskip
{\bf Claim}\qua {\sl The subspace of $H_1(\widetilde{T}, \partial \widetilde{T})$
 on which $f$ acts trivially has dimension at least 4.}

\proof
By Property (1) above, it is enough to show that
 the vectors $[\delta_1], [\delta_1^*],$ $[\delta_2], [\delta_2^*]$
 are linearly independent in $H_1(\widetilde{T}, \partial \widetilde{T})$.
 Let $V_i$ be the space generated by $[\delta_i]$ and $[\delta_i^*]$.
 It follows from Property (2) that $\dim(V_i) = 2$.
 By Property (3), we have $I(v_1, v_2) = 0$ for any $v_1 \in V_1$
 and $v_2 \in V_2$.
 The intersection form $I$ restricted to $V_2$ is a non-zero multiple
 of  the form
  $ \left( \begin{array}{ll} 0 & 1\\  -1 & 0 \end{array} \right) $,
 which is non-singular.
 So, for any $v_2 \in V_2$, there is an element $v_2^* \in V_2$
 such that $I(v_2, v_2^*) \neq 0$.
 Therefore $V_1 \cap V_2 = \emptyset$, so
 the four vectors are linearly independent, and the claim follows.
\endproof

 Each lift of the puncture $\alpha_1$ unwraps twice in $\widetilde{T_1}$
 and once in $\widetilde{T_2}$.
 Therefore each lift of $\alpha_1$ unwraps twice
 in $\widetilde{T}$; similarly, each lift of $\alpha_2$ unwraps
 twice in $\widetilde{T}$.
 Hence there is an induced manifold cover
 $\widetilde{T(2)}$ of $T(2)$ obtained by filling
 in the punctures of $\widetilde{T}$.
 There is then an induced manifold cover $\widetilde{M_f}$ of $M_f$,
 and by Proposition \ref{b1orb},
 $b_1(\widetilde{M_f}) \geq 4 + 1$.

 The proof of the general result is similar. We start
 with an arbitrary positive integer $k$, and
 replace $T$ with a $k$--times punctured torus.

 We then obtain, for each $i\leq k$, a cover $\widetilde{T}_i$ of $T$,
 such that each puncture of $T$ unwraps once or twice in $\widetilde{T}_i$.
 We construct fixed classes
 $[\delta_i], [\delta_i^*] \in H_1(\widetilde{T}_i, \partial \widetilde{T}_i)$
 with algebraic intersection number 2,
 so that the projection of $\delta_i \cup \delta_i^*$ to $T$ is disjoint from
 the projection of $\delta_j \cup \delta_j^*$ whenever $i \neq j$ (see Figure 4).

\begin{figure}[ht!]
\cl{\epsfxsize\hsize\epsfbox{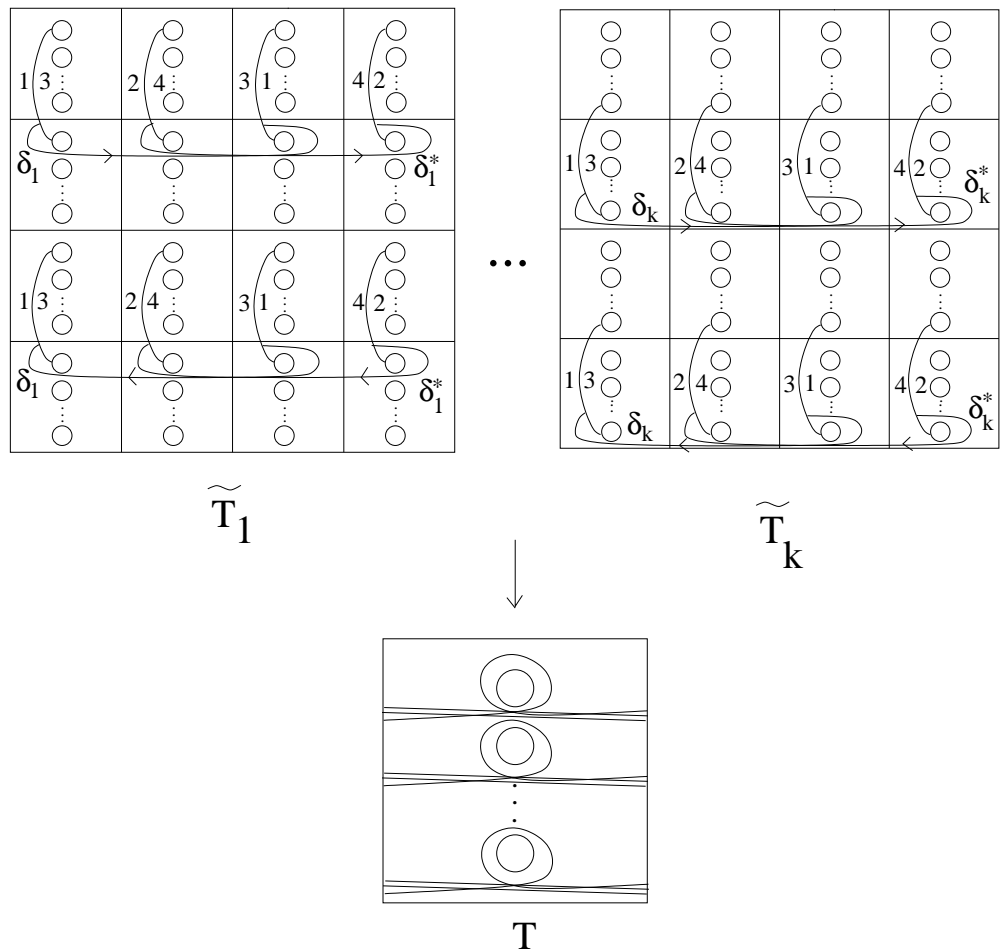}}
\caption{Each boundary component of $T$ gives rise to a different cover.}
\label{coverk}
\end{figure}

 By an argument similar to the one given in the
 $k= 2$ case, we conclude that there is a $2k$--dimensional space
 in $H_1(\widetilde{T}, \partial \widetilde{T})$
 on which $f$ acts trivially.
 Since every puncture of $T$ unwraps twice in $\widetilde{T}$,
 there is an induced manifold cover
 $\widetilde{T(2)}$ of $T(2)$ obtained by filling
 in the punctures of $\widetilde{T}$.
 Therefore there is an induced bundle cover of $M_f$,
 and by Proposition \ref{b1orb},
 $b_1(\widetilde{M_f}) \geq 2k + 1$.
 Since $k$ is an arbitrary positive integer, the result follows.
\end{proof}

\begin{proof}[Proof of Theorem \ref{main}]
This is an application of Lemma \ref{vb1=infty} and Corollary
 \ref{fibundlehomo}
\end{proof}

\section{Proof of Theorem \ref{nonhyper}}

We sketch the proof that $M$ has a finite cover $\widetilde{M}$
 with $b_1(\widetilde{M}) > b_1(M)$.  The generalization to $vb_1(M) = \infty$
 then follows by direct analogy with the proof of Theorem \ref{main}.
 Recall the construction of the cover $\widetilde{F}$ of
 $F$ in the case of a hyper-elliptic monodromy:
  we remove a neighborhood of the fixed points of $\tau$
 to obtain a punctured surface $F^-$. 
 The surface $F^-$ double covers a planar surface $P$;
 we construct a punctured torus $T$
 which double covers $P$, and then a 16--fold cover $\widetilde{T}$ of $T$.
 The cover $\widetilde{F}$ of $F$ corresponds to
 $\pi_1 \widetilde{T} \cap \pi_1 F^-$.
 A loop $\delta  \subset \widetilde{T}$
 is constructed, whose full pre-image
 $\tilde{\delta}$ in $\widetilde{F}$ represents a homology class which 
 is fixed by (a power of) any element of $H= <D_{x_1}, ..., D_{x_{2g}}>$.
 
The covers $T$ and $\widetilde{T}$ of $P$ are not characteristic.
Any element $h$ of $H$ sends $T$ to a cover $hT$ of $P$, and
 $\widetilde{T}$ to a cover $h \widetilde{T}$ of $h T$;
 let $h_0 = id, h_1, ..., h_n \in H$ denote the elements necessary
 for a full orbit of $\widetilde{T}$.
 Let $K_j \subset H_1(h_j \widetilde{T}, \partial h_j \widetilde{T})$
 denote the
 kernel of the projection to $H_1(h_j T, \partial h_j T)$.
 By construction, we have $\delta \in K_0$.
 Let $\gamma$ be the loop pictured in Figure 1, and
 let $p \gamma$ denote the projection of $\gamma$ to $P$.
 We claim that every component of the pre-image of $p \gamma$
 in $h_j \widetilde{T}$ has intersection number 0 with
 every class in $K_j$: this may be checked
 by constructing an explicit basis for the $K_j$'s.

 Now, fix an element $h \in H$.
 The $K_j$'s are permuted by $H$, so $h [\delta]$ 
 has 0 intersection number with each component of the pre-image of $p \gamma$
 in $h \widetilde{T}$.
 Therefore,  every component of the pre-image
 of $\gamma$ in $h \widetilde{F}$ has 0 intersection number
 with $h [\tilde{\delta}]$.
 Since the pre-images of $\gamma$ unwrap at most 8 times,
 we see that $D_{\gamma}^8$ lifts to Dehn twists in $h \widetilde{F}$,
 and fixes $h [\tilde{\delta}]$.
 Therefore, the action of $f$ on $[\tilde{\delta}]$
 is unchanged if we remove all the $D_{\gamma}^8$'s,
 and we deduce from the hyper-elliptic case that, for some integer $m$,
 $f^m$ lifts to a map $\widetilde{f^m}$ such that
 $\widetilde{f^m}_* [\tilde{\delta}] = [\tilde{\delta}]$.
 \endproof

\section{Proof of Theorem \ref{badex}}

 Let $K$ be the knot $9_{32}$ in Rolfsen's tables,
 and let $M = S^3 - K$.  
 The computer program SnapPea shows that $M$ has no symmetries.
 A knot complement is said to have \textit{hidden symmetries}
 if it is an irregular cover of some orbifold.
 In our example, $M$ has no hidden symmetries, since by \cite{NR},
 a hyperbolic knot complement with hidden symmetries must have cusp parameter
 in $\mathbb{Q}(\sqrt{-1})$ or $\mathbb{Q}(\sqrt{-3})$,
 but it is shown in \cite{Ri} that the cusp field of $M$ has degree 29.

 Since $M$ has no symmetries or hidden symmetries,
 and is non-arithmetic (see \cite{R1}), it follows
 from results of Margulis that $M$ is the unique 
 minimal orbifold in its commensurability class.

 Let $M(0,n)$ be the orbifold filling on $K$ obtained by
 setting the $n$-th power of the longitude to the identity. 
 Then, by Corollary 3.3 of \cite{R2}, if $n$ is large enough,
 $M(0,n)$ is a hyperbolic orbifold
 which is minimal in its commensurability class.
 We choose a large $n$ which satisfies this condition and is odd.
 
 Since $K$ has monic Alexander polynomial and fewer than 11 crossings,
 it is fibered (see \cite{K}), and therefore
 $M(0,n)$ is 2--orbifold bundle over $S^1$.
 This orbifold bundle is finitely covered by a manifold which
 fibers over $S^1$;
 let $f\co F \rightarrow F$ denote the monodromy of this fibration.  
 We claim that no power of $f$ lifts to become hyper-elliptic
 in any cover of $F$.
  
 For suppose such a cover $\tilde{F}$ of $F$ exists.
 Then there is an associated cover $\widetilde{M(0,n)}$
 of $M(0,n)$, and an involution $\tau$ on $\widetilde{M(0,n)}$
 with one-dimensional fixed point set.  The quotient
 $Q = \widetilde{M(0,n)}/\tau$ is an orbifold whose singular
 set is a link labeled 2, which is commensurable with $M(0,n)$.
 By minimality, $Q$ must cover $M(0,n)$.  But this is impossible,
 since every torsion element of $M(0,n)$ has odd order,
 by our choice of $n$.
\endproof

\section{Appendix:  Constructing Covers of Punctured Tori}
 We review here the relevant material from \cite{M}.  This builds on work of
 Baker (\cite{B1}, \cite{B2}).

 We are given a punctured torus $T$ and a monodromy $f$,
 and we wish to find finite covers of $T$ to which $f$ lifts.
 Let $x$ and $y$ be the generators for $\pi_1 T$ pictured in
 Figure 5.
 Let $r$ and $s$ be positive integers, and let
 $\hat{T}$ be the  $rs$--fold cover of $T$ associated to the kernel
 of the map $\phi \co  \pi_1(T) \rightarrow \mathbb{Z}_r \times \mathbb{Z}_s$,
 with $\phi([x]) = (1,0)$ and $\phi([y]) = (0,1)$ (see Figure 5).

\begin{figure}[ht!]
\cl{\epsfbox{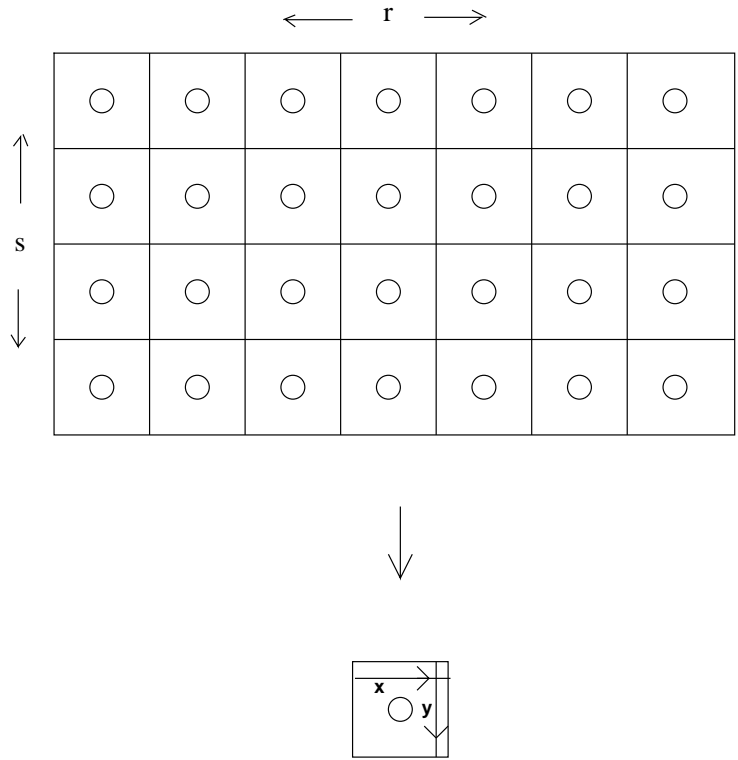}}
\caption{The cover $\hat{T}$ of $T$}
\end{figure}

 Now we create a new cover, $\tilde{T}$, of $T$ by making vertical cuts in
 each row of $\hat{T}$,
 and gluing the left side of each cut to the right side of another cut in
 the same row.  An example is pictured in Figure 6,
 where the numbers
 in each row indicate how the edges are glued.

\begin{figure}[ht!]
\cl{\epsfxsize.9\hsize\epsfbox{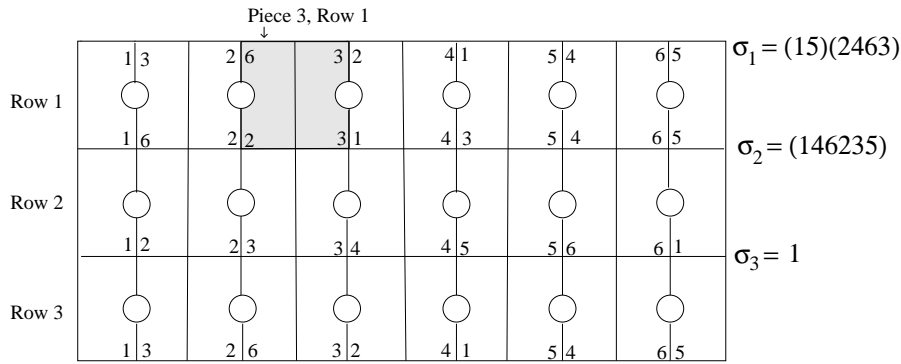}}
\caption{The permutations encode the combinatorics of the gluing}
\label{comb}
\end{figure}

We now introduce some notation to describe the cuts of $\tilde{T}$
 (see Figure 6). 
 $\tilde{T}$ is naturally divided into rows, which we label
 $1, ..., s$.  The cuts divide each row
 into pieces, each of which is a square minus two half-disks;
 we number them $1, ..., r$.  If we slide each point in the top half of
 the $i^{th}$ row through the cut to its right, we induce a permutation on
 the pieces $\lbrace 1, ..., r \rbrace$, which we denote $\sigma_i$.
  Thus the cuts on $\tilde{T}$
 may be encoded by elements $\sigma_1, ..., \sigma_s \in S_r$, the
 permutation group on $r$ letters.

 Let $D_x$ and $D_y$ be the right-handed Dehn twists
 in $x$ and $y$, which generate the mapping class group of $T$.
 We observe that, regardless of the choice of $\sigma_i$'s,
 $D_y^s$ lifts to a product of Dehn twists in $\tilde{T}$.
 It will be useful to have a condition on the $\sigma_i$'s which
 will guarantee that $D_x$ lifts to $\tilde{T}$.
 The following lemma (in slightly different form)
 appears in \cite{M}.

 \begin{lem} \label{lift}
$D_x$ lifts to $\tilde{D}_x\co  \tilde{T} \rightarrow \tilde{T}$ if\\
{\rm(1)}\qua $\sigma_1 ... \sigma_i$
 commutes with $\sigma_{i+1}$ for $i = 1, ..., s-1$, and\\
{\rm(2)}\qua $\sigma_1 ... \sigma_s = 1$.\\
Moreover, if these conditions are satisfied,
 then we may choose $\tilde{D}_x$ so that its
action on the interior of the $i$th row
 of $\tilde{T}$ corresponds to the permutation
 $\sigma_1 ... \sigma_i$.
\end{lem}

\begin{proof}
  We shall attempt to lift $D_x$ explicitly
 to a sequence of ``fractional Dehn twists'' along the rows of $\tilde{T}$.
  Let $\tilde{x}_i$ denote the disjoint union of the lifts of $x$
 to the $i^{th}$ row of $\tilde{T}$.
  We first attempt to lift $D_x$ to row 1, twisting one slot to the right
 along $\tilde{x}_1$.
 Considering the effect of this
 action on the bottom half of row 1, we find the cuts there
are now matched up according to the permutation
 $\sigma_1^{-1}\sigma_2\sigma_1$.  Thus, for $D_x$ to lift to row 1
 we assume $\sigma_1$ and $\sigma_2$ commute.  We now twist along
 $\tilde{x}_2$.  The top halves of the squares in row 2 are moved
 according to the permutation
 $\sigma_1\sigma_2$, and the lift will extend to all of row 2 if
 $\sigma_3$ commutes with $\sigma_1\sigma_2$.  We continue in this manner,
 obtaining the conditions in 1.  After we twist through
 $\tilde{x}_n$, we need to be back where we started in row 1;
 if the permutations satisfy
 the additional condition $\sigma_1\sigma_2...\sigma_s = 1$,
 then this is the case, and we have succeeded in lifting $D_x$.
 Note that in the course of constructing the lift, we
 have also verified the last assertion of the lemma.
\end{proof}

For the purposes of this paper,
we restrict attention to the case $s = 4$.
Consider the subgroup $J = <D_x, D_y^4>$ of the mapping
 class group of $T$.
If $\sigma_1, ..., \sigma_4$ satisfy the conditions of Lemma
 \ref{lift}, then
 any element of $J$ lifts to $\tilde{T}$.   What makes this useful is the
 following lemma.

\begin{lem}\label{fi}
The subgroup $J$ has finite index in the mapping class group of $T$.
\end{lem}

\begin{proof}
The mapping class group of $T$ may be indentified with $SL_2(\mathbb{Z})$,
 and under this indentification, $J$ is the group generated
 by  $ \left( \begin{array}{ll} 1 & 0\\  4 & 1 \end{array} \right) $
 and $ \left( \begin{array}{ll} 1 & 1\\  0 & 1 \end{array} \right) $.

Let $\gamma = \left( \begin{array}{ll} \sqrt{2} & 0\\  0 & \frac{1}{\sqrt{2}} \end{array} \right)$.  Then $\gamma$ conjugates the generators of $J$
 to  $ \left( \begin{array}{ll} 1 & 0\\  2 & 1 \end{array} \right) $
 and $ \left( \begin{array}{ll} 1 & 2\\  0 & 1 \end{array} \right) $,
 which are well known to generate the kernel of the reduction
 map from $SL_2(\mathbb{Z})$ to $SL_2(\mathbb{Z}/2)$.
 Therefore $J$ is a finite co-area lattice in $SL_2(\mathbb{R})$,
 and therefore it has finite index in $SL_2(\mathbb{Z})$.
\end{proof}

The next lemma shows that with some additional hypotheses on the
 $\sigma_i$'s we are also guaranteed that the lifts of elements
 of $J$ fix non-peripheral homology classes of $\tilde{T}$.

\begin{lem} \label{int0}
Let $\tilde{T}$ be as constructed above, and suppose $\sigma_2 = \sigma_1^{-1}$
 and $\sigma_4 = \sigma_3^{-1}$.  Let $f$ be an element of $J$.  Then\\
{\rm(i)}\qua  $f$ lifts to an automorphism
 $\tilde{f}\co  \tilde{T} \rightarrow \tilde{T}$, and\\
{\rm(ii)}\qua  For every non-peripheral loop $\ell$ in Row 2, there is a loop $\ell^*$ in Row 4, such that $\tilde{f}_* [\ell \cup \ell^*] = [\ell \cup \ell^*] \neq [0] \in H_1(\tilde{T}, \partial \tilde{T})$.
\end{lem}

\begin{proof}
 Assertion (i) is an immediate consequence of Lemma \ref{lift}. 
 To prove Assertion (ii), we explicitly construct
 the loop $\ell^*$, so that it intersects the same components
 of $\tilde{y}$ as $\ell$ does, but with opposite orientations.
 Figure 7 indicates the procedure for doing this.

 Therefore $[\ell \cup \ell^*]$ has 0 intersection number with
 each component of $\tilde{y}$,
 and so it is fixed homologically by $\tilde{D}^4_y$.  Moreover,
 $\ell \cup \ell^*$ is entirely contained in Rows 2 and 4,
 and Lemma \ref{lift} implies that the action of $\tilde{D}_x$
 is trivial there, so $[\ell \cup \ell^*]$ is also fixed by $\tilde{D}_x$,
 and by every element of $J$.
\end{proof}

\newpage
\begin{figure}[ht!] \label{cancel}
\cl{\epsfxsize4.8in\epsfbox{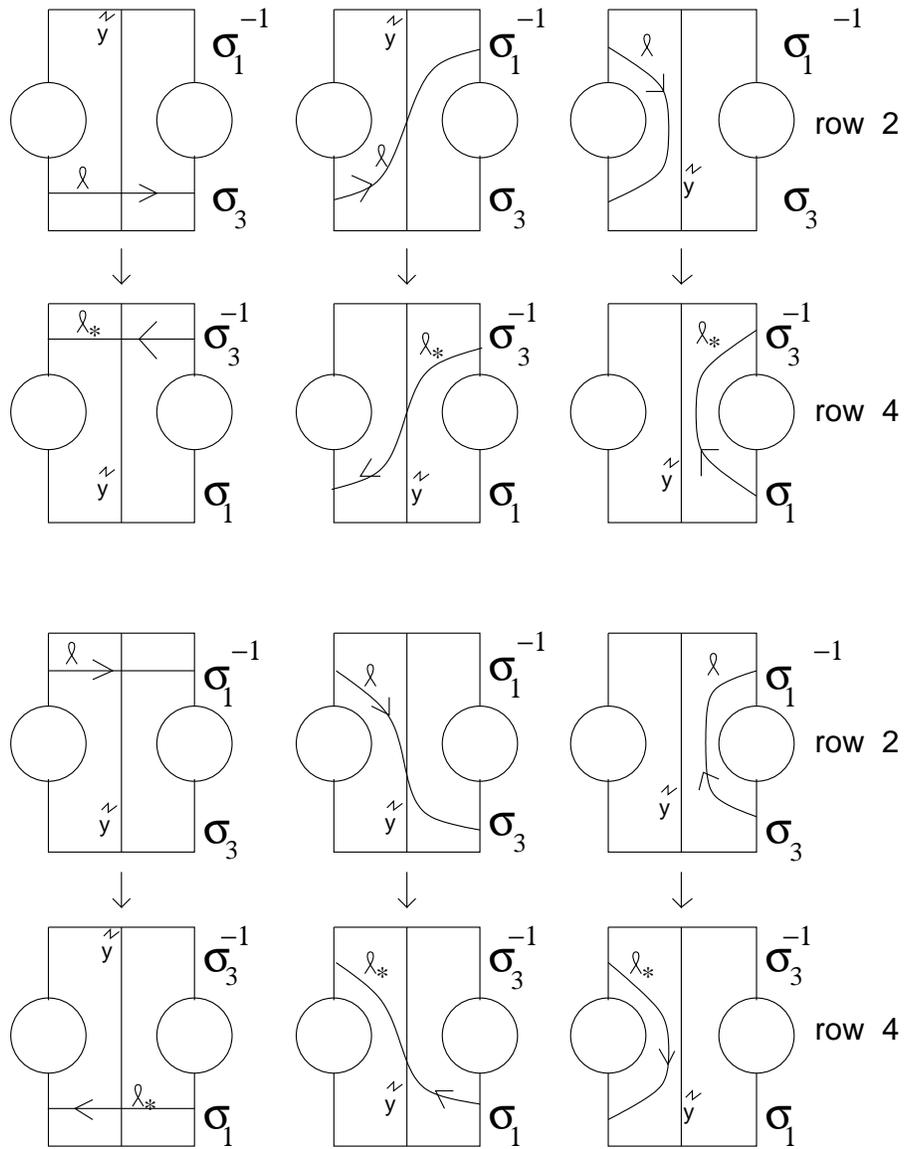}}
\caption{Corresponding to each segment of $\ell$,
 we construct a corresponding segment of $\ell_*$.}
\end{figure}


\newpage

\begin{thebibliography}

\bibitem{B1} {\bf M Baker}, {\it Covers of Dehn fillings on
 once-punctured torus bundles}, {Proc. Amer. Math. Soc.} {105} (1989)
 747--754

\bibitem{B2} {\bf M Baker}, {\it Covers of Dehn fillings on
 once-punctured torus bundles II}, {Proc. Amer. Math. Soc.} {110}
 (1990) 1099--1108

\bibitem{G} {\bf D Gabai}, {\it On 3--manifolds finitely covered by
  surface bundles}, from: ``Low-dimensional Topology and Kleinian Groups
  (Coventry/Durham, 1984)'', LMS Lecture Note Series 112, Cambridge
  University Press (1986)

\bibitem{H} {\bf S\,P Humphries}, {\it Generators for the mapping
 class group}, from: ``Topology of Low-Dimensional Manifolds'', Proceedings
 of the Second Sussex Conference, 1977, Lecture Notes in Mathematics
 722, Springer--Verlag, Berlin (1979)

\bibitem{K} {\bf T Kanenobu}, {\it The augmentation subgroup of a
 pretzel link}, Mathematics Seminar Notes, Kobe University, {7} (1979)
 363--384

\bibitem{M} {\bf J\,D Masters}, {\it Virtual homology of surgered torus
bundles"}, to appear in Pacific J. Math.

\bibitem{NR} {\bf W\,D Neumann}, {\bf A\,W Reid}, {\it Arithmetic of hyperbolic
 3--manifolds}, Topology '90, de Gruyter (1992) 273--309

\bibitem{O} {\bf J-P Otal}, {\it Thurston's hyperbolization of Haken
manifolds}, from: ``Surveys in Differential Geometry, Vol. III'',
(Cambridge MA 1996), Int. Press, Boston MA (1998) 77--194

\bibitem{R1} {\bf A\,W Reid}, {\it Arithmeticity of knot complements},
 {J. London Math. Soc.} {43} (1991) 171--184

\bibitem{R2} {\bf A\,W Reid}, {\it Isospectrality and commensurability
 of arithmetic hyperbolic 2-- and 3--manifolds}, {Duke Math. J.}
 {65} (1992), no. 2, 215--228

\bibitem{Ri} {\bf R Riley}, {\it Parabolic representations and
 symmetries of the knot $9_{32}$}, from: ``Computers and Geometry and
 Topology'', (M\,C Tangora, editor), Lecture Notes in Pure and Applied
 Math. {114}, Dekker (1988) 297--313

\bibitem{T} {\bf W Thurston}, {\it A norm for the homology of 3--manifolds},
 {Mem. Amer. Math. Soc.} {59} (1986)
 99--130

\bibitem{W} {\bf F Waldhausen}, {\it On irreducible 3--manifolds which are
 sufficiently large}, Ann. of Math. {87} (1968)
 195--203

\end{thebibliography}
\end{document}